\documentclass[reqno]{amsart}
\usepackage{url}
\usepackage[dvips]{graphicx}

\newtheorem{thm}{Theorem}
\newtheorem{lem}{Lemma}
\newtheorem{definition}{Definition}
\begin{document}
\bibliographystyle{plain}

\title[Constructive proof of Brouwer's fixed point theorem]{Constructive proof of Brouwer's fixed point theorem for sequentially locally non-constant functions}

\author{Yasuhito Tanaka}
\address{Faculty of Economics, Doshisha University, Kamigyo-ku, Kyoto, 602-8580, Japan}
\email{yasuhito@mail.doshisha.ac.jp}
\thanks{This research was partially supported by the Ministry of Education, Science, Sports and Culture of Japan, Grant-in-Aid for Scientific Research (C), 20530165.}\address{Faculty of Economics, Doshisha University, Kamigyo-ku, Kyoto, 602-8580, Japan}

\date{}



\begin{abstract}  
We present a constructive proof of Brouwer's fixed point theorem for uniformly continuous and sequentially locally non-constant functions based on the existence of approximate fixed points. And we will show that Brouwer's fixed point theorem for uniformly continuous and sequentially locally non-constant functions implies Sperner's lemma for a simplex. Since the existence of approximate fixed points is derived from Sperner's lemma, our Brouwer's fixed point theorem is equivalent to Sperner's lemma. 
\end{abstract}

\maketitle

\section{Introduction}

It is well known that Brouwer's fixed point theorem can not be constructively proved\footnote{\cite{kel} provided a \emph{constructive} proof of Brouwer's fixed point theorem. But it is not constructive from the view point of constructive mathematics \'{a} la Bishop. It is sufficient to say that one dimensional case of Brouwer's fixed point theorem, that is, the intermediate value theorem is non-constructive. See \cite{br} or \cite{da}.}. Sperner's lemma which is used to prove Brouwer's theorem, however, can be constructively proved. Some authors have presented an approximate version of Brouwer's theorem using Sperner's lemma. See \cite{da} and \cite{veld}.  Thus, Brouwer's fixed point theorem is constructively, in the sense of constructive mathematics \'{a} la Bishop, proved in its approximate version.

Also \cite{da} states a conjecture that a uniformly continuous function $f$ from a simplex into itself, with property that each open set contains a point $x$ such that $x\neq f(x)$, which means $|x-f(x)|>0$, and also at every point $x$ on the boundaries of the simplex $x\neq f(x)$, has an exact fixed point. We call such a property of functions \emph{local non-constancy}. In this note we present a partial answer to Dalen's conjecture.

 Recently \cite{berger} showed that the following theorem is equivalent to Brouwer's fan theorem.
\begin{quote}
Each uniformly continuous function $f$ from a compact metric space $X$ into itself with at most one fixed point and approximate fixed points has a fixed point.
\end{quote}
\clearpage
By reference to the notion of \emph{sequentially at most one maximum} in \cite{berg} we require a condition that a function $f$ is \emph{sequentially locally non-constant}, and will show the following result.
\begin{quote}
Each uniformly continuous function $f$ from an $n$-dimensional simplex into itself, which is sequentially locally non-constant and has approximate fixed points, has a fixed point,
\end{quote}
\noindent without the fan theorem. Sequential local non-constancy is stronger than the condition in \cite{da} (local non-constancy), and is different from the condition that a function has \emph{at most one fixed point} in \cite{berg}.

\cite{orevkov} constructed a computably coded continuous function $f$ from the unit square into itself, which is defined at each computable point of the square, such that $f$ has no computable fixed point. His map consists of a retract of the computable elements of the square to its boundary followed by a rotation of the boundary of the square. As pointed out by \cite{hirst}, since there is no retract of the square to its boundary, his map does not have a total extension.

In the next section we present our theorem and its proof. In Section 3 we will derive Sperner's lemma from Brouwer's fixed point theorem for uniformly continuous and sequentially locally non-constant functions.

\section{Brouwer's fixed point theorem for sequentially locally non-constant functions}

We consider an $n$-dimensional simplex $\Delta$ as a compact metric space. Let $x$ be a point in $\Delta$, and consider a uniformly continuous function $f$ from $\Delta$ into itself.

 According to \cite{da} and \cite{veld} $f$ has an approximate fixed point. It means 
\[\mathrm{For\ each}\ \varepsilon>0\ \mathrm{there\ exists}\ x\in \Delta\ \mathrm{such\ that}\ |x-f(x)|<\varepsilon.\]
Since $\varepsilon>0$ is arbitrary,
\[\inf_{x\in \Delta}|x-f(x)|=0.\]

The definition of local non-constancy of functions in a simplex $\Delta$ is as follows;
\begin{definition}(Local non-constancy of function)
\begin{enumerate}
	\item At a point $x$ on the faces (boundaries) of $\Delta$ $f(x)\neq x$. This means $f_i(x)>x_i$ or $f_i(x)<x_i$ for at least one $i$, where $x_i$ and $f_i(x)$ denote the $i$-th components of $x$ and $f(x)$.
	\item And in any open set in $\Delta$ there exists a point $x$ such that $f(x)\neq x$.
\end{enumerate}
\end{definition}

On the other hand the notion that $f$ has at most one fixed point by \cite{berger} is defined as follows;
\begin{definition}[At most one fixed point]
For all $x, y\in \Delta$, if $x\neq y$, then $f(x)\neq x$ or $f(y)\neq y$.
\end{definition}

By reference to the notion of \emph{sequentially at most one maximum} in \cite{berg},  we define the property of \emph{sequential local non-constancy}.

First we recapitulate the compactness (total boundedness with completeness) of a set in constructive mathematics. $\Delta$ is totally bounded in the sense that for each $\varepsilon>0$ there exists a finitely enumerable $\varepsilon$-approximation to $\Delta$\footnote{A set $S$ is finitely enumerable if there exist a natural number $N$ and a mapping of the set $\{1, 2, \dots, N\}$ onto $S$.}. An $\varepsilon$-approximation to $\Delta$ is a subset of $\Delta$ such that for each $x\in \Delta$ there exists $y$ in that $\varepsilon$-approximation with $|x-y|<\varepsilon$. According to Corollary 2.2.12 of \cite{bv} we have the following result.
\begin{lem}
For each $\varepsilon>0$ there exist totally bounded sets $H_1, H_2, \dots, H_n$, each of diameter less than or equal to $\varepsilon$, such that $\Delta=\cup_{i=1}^nH_i$.\label{closed}
\end{lem}
Since $\inf_{x\in \Delta}|x-f(x)|=0$, we have $\inf_{x\in H_i}|x-f(x)|=0$ for some $H_i\subset \Delta$.

The definition of sequential local non-constancy is as follows;
\begin{definition}(Sequential local non-constancy of functions)
There exists $\bar{\varepsilon}>0$ with the following property. For each $\varepsilon>0$ less than or equal to $\bar{\varepsilon}$ there exist totally bounded sets $H_1, H_2, \dots, H_m$, each of diameter less than or equal to $\varepsilon$, such that $\Delta=\cup_{i=1}^m H_i$, and if for all sequences $(x_n)_{n\geq 1}$, $(y_n)_{n\geq 1}$ in each $H_i$, $|f(x_n)-x_n|\longrightarrow 0$ and $|f(y_n)-y_n|\longrightarrow 0$, then $|x_n-y_n|\longrightarrow 0$.
\end{definition}

Now we show the following lemma, which is based on Lemma 2 of \cite{berg}.
\begin{lem}
Let $f$ be a uniformly continuous function from $\Delta$ into itself. Assume $\inf_{x\in H_i}f(x)=0$ for some $H_i\subset \Delta$ defined above. If the following property holds:
\begin{quote}
For each $\delta>0$ there exists $\varepsilon>0$ such that if $x, y\in H_i$, $|f(x)-x|<\varepsilon$ and $|f(y)-y|<\varepsilon$, then $|x-y|\leq \delta$.
\end{quote}
Then, there exists a point $z\in H_i$ such that $f(z)=z$, that is, $f$ has a fixed point. \label{fix0}
\end{lem}
\begin{proof}
Choose a sequence $(x_n)_{n\geq 1}$ in $H_i$ such that $|f(x_n)-x_n|\longrightarrow 0$. Compute $N$ such that $|f(x_n)-x_n|<\varepsilon$ for all $n\geq N$. Then, for $m, n\geq N$ we have $|x_m-x_n|\leq \delta$. Since $\delta>0$ is arbitrary, $(x_n)_{n\geq 1}$ is a Cauchy sequence in $H_i$, and converges to a limit $z\in H_i$. The continuity of $f$ yields $|f(z)-z|=0$, that is, $f(z)=z$.
\end{proof}

Next we show the following theorem, which is based on Proposition 3 of \cite{berg}.
\begin{thm}
Each uniformly continuous function $f$ from an $n$-dimensional simplex into itself, which is sequentially locally non-constant and has approximate fixed points, has a fixed point.
\end{thm}
\begin{proof}
Choose a sequence $(z_n)_{n\geq 1}$ in $H_i$ defined above such that $|f(z_n)-z_n|\longrightarrow 0$. In view of Lemma \ref{fix0} it is enough to prove that the following condition holds.
\begin{quote}
For each $\delta>0$ there exists $\varepsilon>0$ such that if $x, y\in H_i$, $|f(x)-x|<\varepsilon$ and $|f(y)-y|<\varepsilon$, then $|x-y|\leq \delta$.
\end{quote}
Assume that the set
\[K=\{(x,y)\in H_i\times H_i:\ |x-y|\geq \delta\}\]
is nonempty and compact\footnote{See Theorem 2.2.13 of \cite{bv}.}. Since the mapping $(x,y)\longrightarrow \max(|f(x)-x|,|f(y)-y|)$ is uniformly continuous, we can construct an increasing binary sequence $(\lambda_n)_{n\geq 1}$ such that
\[\lambda_n=0\Rightarrow \inf_{(x,y)\in K}\max(|f(x)-x|,|f(y)-y|)<2^{-n},\]
\[\lambda_n=1\Rightarrow \inf_{(x,y)\in K}\max(|f(x)-x|,|f(y)-y|)>2^{-n-1}.\]
It suffices to find $n$ such that $\lambda_n=1$. In that case, if $|f(x)-x|<2^{-n-1}$, $|f(y)-y|<2^{-n-1}$, we have $(x,y)\notin K$ and $|x-y|\leq \delta$. Assume $\lambda_1=0$. If $\lambda_n=0$, choose $(x_n, y_n)\in K$ such that $\max(|f(x_n)-x_n|, |f(y_n)-y_n|)<2^{-n}$, and if $\lambda_n=1$, set $x_n=y_n=z_n$. Then, $|f(x_n)-x_n|\longrightarrow 0$ and $|f(y_n)-y_n|\longrightarrow 0$, so $|x_n-y_n|\longrightarrow 0$. Computing $N$ such that $|x_N-y_N|<\delta$, we must have $\lambda_N=1$. We have completed the proof.
\end{proof}

\section{From Brouwer's fixed point theorem for sequentially locally non-constant functions to Sperner's lemma}

\begin{figure}[t]
\begin{center}
\includegraphics[height=7.5cm]{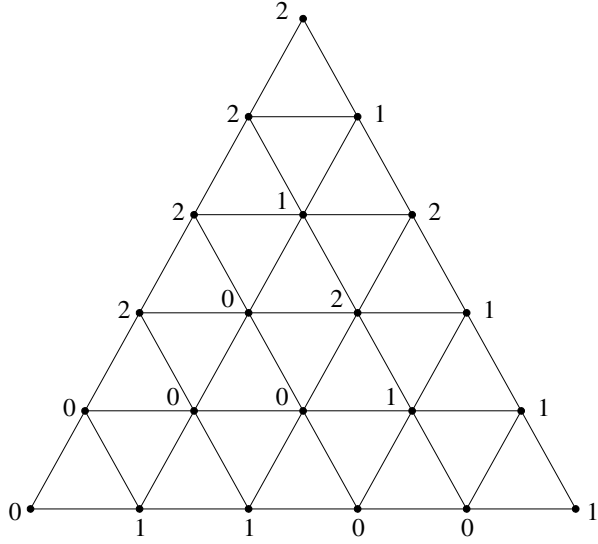}
\end{center}
	\vspace*{-.3cm}
	\caption{Partition and labeling of 2-dimensional simplex}
	\label{tria2}
\end{figure}

In this section we will derive Sperner's lemma from Brouwer's fixed point theorem for uniformly continuous and sequentially locally non-constant functions. Let $\Delta$ be an $n$-dimensional simplex. Denote a point on $\Delta$ by $x$. Consider a function $f$ from $\Delta$ into itself. Partition $\Delta$ in the way depicted in Figure \ref{tria2} for a 2-dimensional simplex. In a 2-dimensional case we divide each side of $\Delta$ in $m$ equal segments, and draw the lines parallel to the sides of $\Delta$. Then, the 2-dimensional simplex is partitioned into $m^2$ triangles. We consider partition of $\Delta$ inductively for cases of higher dimension. In a 3 dimensional case each face of $\Delta$ is a 2-dimensional simplex, and so it is partitioned into $m^2$ triangles in the above mentioned way, and draw the planes parallel to the faces of $\Delta$. Then, the 3-dimensional simplex is partitioned into $m^3$ trigonal pyramids. And similarly for cases of higher dimension. Let $K$ denote the set of small $n$-dimensional simplices of $\Delta$ constructed by partition. Vertices of these small simplices of $K$ are labeled with the numbers 0, 1, 2, $\dots$, $n$ subject to the following rules. 
\begin{enumerate}
\item The vertices of $\Delta$ are respectively labeled with 0 to $n$. We label a point $(1,0, \dots, 0)$ with 0, a point $(0,1,0, \dots, 0)$ with 1, a point $(0,0,1 \dots, 0)$ with 2, $\dots$, a point $(0,\dots, 0,1)$ with $n$. That is, a vertex whose $k$-th coordinate ($k=0, 1, \dots, n$) is $1$ and all other coordinates are 0 is labeled with $k$ for all $k\in \{0, 1, \dots, n\}$. 

\item If a vertex of a simplex of $K$ is contained in an $n-1$-dimensional face of $\Delta$, then this vertex is labeled with some number which is the same as the number of a vertex of that face.

\item If a vertex of a simplex of $K$ is contained in an $n-2$-dimensional face of $\Delta$, then this vertex is labeled with some number which is the same as the number of a vertex of that face. And similarly for cases of lower dimension.

\item A vertex contained inside of $\Delta$ is labeled with an arbitrary number among 0, 1, $\dots$, $n$.
\end{enumerate}

Denote the vertices of an $n$-dimensional simplex of $K$ by $x^0, x^1, \dots, x^n$, the $j$-th coordinate of $x^i$ by $x^i_j$, and denote the label of $x^i$ by $l(x^i)$. Let $\tau$ be a positive number which is smaller than $x^i_{l(x^i)}$ for all $x^i$, and define a function $f(x^i)$ as follows\footnote{We refer to \cite{yoseloff} about the definition of this function.}:
\[f(x^i)=(f_0(x^i), f_1(x^i), \dots, f_n(x^i)),\]
and
\begin{equation}
f_j(x^i)=\left\{
\begin{array}{ll}
x^i_j-\tau&\mathrm{for}\ j=l(x^i),\\
x^i_j+\frac{\tau}{n}&\mathrm{for}\ j\neq l(x^i).\label{e0}
\end{array}
\right.
\end{equation}
$f_j$ denotes the $j$-th component of $f$. From the labeling rules we have $x^i_{l(x^i)}>0$ for all $x^i$, and so $\tau>0$ is well defined. Since $\sum_{j=0}^nf_j(x^i)=\sum_{j=0}^nx^i_j=1$, we have
\[f(x^i)\in \Delta.\]
We extend $f$ to all points in the simplex by convex combinations on the vertices of the simplex. Let $z$ be a point in the $n$-dimensional simplex of $K$ whose vertices are $x^0, x^1, \dots, x^n$. Then, $z$ and $f(z)$ are expressed as follows;
\[z=\sum_{i=0}^n\lambda_ix^i,\ \mathrm{and}\ f(z)=\sum_{i=0}^n\lambda_if(x^i),\ \lambda_i\geq 0,\ \sum_{i=0}^n\lambda_i=1.\]

It is clear that $f$ is uniformly continuous. We verify that $f$ is sequentially locally non-constant.
\begin{enumerate}
\item Let $z$ be a point in an $n$-dimensional simplex $\delta^n$. Assume that no vertex of $\delta^n$ is labeled with $i$. Then
\begin{equation}
f_i(z)=\sum_{j=0}^n\lambda_jf_i(x^j)=z_i+\left(1+\frac{1}{n}\right)\tau.\label{e8}
\end{equation}
Then, there exists no sequence $(z_m)_{m\geq 1}$ such that $|f(z_m)-z_m|\longrightarrow 0$ in $\delta^n$.
\item Assume that $z$ is contained in a fully labeled $n$-dimensional simplex $\delta^n$, and rename vertices of $\delta^n$ so that a vertex $x^i$ is labeled with $i$ for each $i$. Then,
\begin{align*}
f_i(z)=\sum_{j=0}^n\lambda_jf_i(x^j)=\sum_{j=0}^n\lambda_jx_i^j+\sum_{j\neq i}\lambda_j\frac{\tau}{n}-\lambda_i\tau=z_i+\left(\frac{1}{n}\sum_{j\neq i}\lambda_j-\lambda_i\right)\tau\ \mathrm{for\ each}\ i.
\end{align*}
Consider sequences $(z_m)_{m\geq 1}$, $(z'_m)_{m\geq 1}$ such that $|f(z_m)-z_m|\longrightarrow 0$ and $|f(z'_m)-z'_m|\longrightarrow 0$.

Let $z_m=\sum_{i=0}^n\lambda(m)_ix^i$ with $\lambda(m)_i\geq 0,\ \sum_{i=0}^n\lambda(m)_i=1$ and $z'_m=\sum_{i=0}^n\lambda'(m)_ix^i$  with $\lambda'(m)_i\geq 0,\ \sum_{i=0}^n\lambda'(m)_i=1$. Then, we have
\[\frac{1}{n}\sum_{j\neq i}\lambda(m)_j-\lambda(m)_i\longrightarrow 0,\ \mathrm{and}\ \frac{1}{n}\sum_{j\neq i}\lambda'(m)_j-\lambda'(m)_i\longrightarrow 0\ \mathrm{for\ all}\ i.\]
Therefore, we obtain
\[\lambda(m)_i\longrightarrow \frac{1}{n+1}, \mathrm{and}\ \lambda'(m)_i\longrightarrow \frac{1}{n+1}.\]
These mean
\[|z_m-z'_m|\longrightarrow 0.\]
\end{enumerate}
Thus, $f$ is sequentially locally non-constant, and it has a fixed point. Let $z^*$ be a fixed point of $f$. We have
\begin{equation}
z^*_i=f_i(z^*)\ \mathrm{for\ all}\ i.\label{e15}
\end{equation}
Suppose that $z^*$ is contained in a small $n$-dimensional simplex $\delta^*$. Let $x^0, x^1, \dots, x^n$ be the vertices of $\delta^*$. Then, $z^*$ and $f(z^*)$ are expressed as
\[z^*=\sum_{i=0}^n\lambda_ix^i\ \mathrm{and}\ f(z^*)=\sum_{i=0}^n\lambda_if(x^i),\ \lambda_i\geq 0,\ \sum_{i=0}^n\lambda_i=1.\]
(\ref{e0}) implies that if only one $x^k$ among $x^0, x^1, \dots, x^n$ is labeled with $i$, we have
\[f_i(z^*)=\sum_{j=0}^n\lambda_jf_i(x^j)=\sum_{j=0}^n\lambda_jx_i^j+\sum_{j\neq k}^n\lambda_j\frac{\tau}{n}-\lambda_k\tau=z_i^*\ \mathrm{(}z_i^*\mathrm{\ is\ the}\ i\mathrm{-th\ coordinate\ of}\ z^*\mathrm{)}.\]
This means
\[\frac{1}{n}\sum_{j\neq k}^n\lambda_j-\lambda_k=0.\]
Then, (\ref{e15}) is satisfied with $\lambda_k=\frac{1}{n+1}$ for all $k$. If no $x^j$ is labeled with $i$, we have (\ref{e8}) with $z=z^*$ and then (\ref{e15}) can not be satisfied. Thus, one and only one $x^j$ must be labeled with $i$ for each $i$. Therefore, $\delta^*$ must be a fully labeled simplex, and so the existence of a fixed point of $f$ implies the existence of a fully labeled simplex.

We have completely proved Sperner's lemma.

\bibliography{yasuhito}

\end{document}